# ON THE PROBABILISTIC RATIONALE OF *I*-DIVERGENCE AND *J*-DIVERGENCE MINIMIZATION

MARIAN GRENDÁR, JR AND MARIÁN GRENDÁR

ABSTRACT. A probabilistic rationale for *I*-divergence minimization (relative entropy maximization), non-parametric likelihood maximization and *J*-divergence minimization (Jeffreys' entropy maximization) criteria is provided.

## 1. INTRODUCTION

Appearance of mathematical expressions known nowadays as *I*-divergence and *J*-divergence dates back to Jeffreys' investigations of invariant priors, [9]. Their statistical and mathematical properties were further explored by [12] and [11], developing an information theoretic approach to statistics. *I*-divergence is since then known also as Kullback-Leibler (KL) information/distance or KL directed divergence.

Despite the important role which *I*-divergence minimization criterion plays in several subjects of mathematical statistics, it lacks rationale, justification. *I*-divergence itself, is usually explicated as being designed to measure how easy it is to tell that data having one distribution do not have some other distribution. Why just *I*-divergence and not some other 'divergence'? In the context of hypothesis testing, Sanov's theorem gives (at least, or at most) an operational meaning to the *I*-divergence. That is probably all what can mathematical statistics say concerning a rationale for *I*-divergence *minimization criterion*.

*J*-divergence, unlike *I*-divergence, has not become popular among statisticians, thus for what it was designed is yet not known.

KL distance, taken with minus sign, is known also as relative entropy or cross-entropy $-\sum_{i=1}^{m} p_i \ln(p_i/q_i)$ of a pmf **p** relative to **q**. With the base measure **q** taken to be uniform pmf, the relative entropy reduces (when maximized) to Shannon's entropy. Its maximization was raised up into a principle status (MaxEnt) by [5], in statistical physics. Later on, Jaynes promoted the MaxEnt principle also in the area of statistical inference, see [8]. Several justification of MaxEnt as a method

1991 *Mathematics Subject Classification*. Primary 60C05; Secondary 60F99, 62A15, 62A99.

Based on Technical Report 2/2000 of Inst. of Measurement Science, Slovak Academy of Sciences, April 2000. Presented at 'Information Theory in Mathematics' conference, Balatonlelle, Hungary, July 4 – 7, 2000. It is a pleasure to thank several participants of the conference, of 'Statistics and Probability Seminar' at FMPH CU Bratislava, of 'MaxEnt 2000' conference for valuable discussions and comments. Special thanks go to Ariel Caticha, Peter Cheeseman, Aleš Gottvald, George Judge, Mike Malioutov, Teddy Seidenfeld, Alberto Solana, Igor Vajda, Benjamin Weiss and Viktor Witkovský for intellectual generosity. Lapses are obviously ours.

for judging under uncertainty were proposed, ranging from a combinatorial argument (Wallis-Jaynes theorem, see [8]), through entropy concentration theorem (see [6]) to several axiomatizations (see for instance [13] or [1]) which single out the Shannon's entropy function but leave unanswered the question recurring in Jaynes' thoughts: 'Just what we are accomplishing when we maximize entropy?', [7]. The combinatorial justifications of MaxEnt are critically assessed in [4].

The article offers a probabilistic rationale for $I$-divergence minimization or equivalently Relative Entropy Maximization (REM) and MaxEnt, as well as $J$-divergence minimization or Jeffreys' Entropy Maximization (JEM). Next section introduces the first elementary setup, asks simple question (Question 1, Q1): *what is the most probable vector of absolute frequencies (occurrences), among vectors allowed by a constraints, which can be generated by a probability distribution (generator)?*, states the main result (Theorem 1) showing that $I$-divergence minimization is an asymptotic form (in the sense of $n \to \infty$) of the Q1. Thus, the only proper place of $I$-divergence minimization is in asymptotic. For a known, finite sample size $n$ it should be replaced by the Q1 method, dubbed *MaxProb* at [4]. The claim of Theorem 1 is illustrated by a numerical calculation (Example 1). Next, it is shown, that another simple probabilistic question (Question 2, Q2): *what is the expected vector of absolute frequencies (occurrences), among vectors allowed by a constraints, which can be generated by a probability distribution (generator)?*, leads asymptotically to $I$-divergence minimization, too. In Section 3 non-parametric likelihood maximization, as a logically inverse problem to those of Section 2, is briefly summarized and illustrated on Example 3. It serves mainly for mediating a move to the rationale for $J$-divergence minimization, presented in Section 4. Section 5 summarizes the main results. Proof of the main Theorem 1 is in Appendix.

## 2. $I$-div minimization - as an asymptotic form of Question 1, 2

Let us consider following general setup.

**Setup 1.** Let $\mathbf{q} = [q_1, q_2, \ldots, q_m]'$ be a pmf, defined on $m$-element support, referred to as *prior generator*.

Let $\mathcal{H}_n$ be a set of all vectors $\{\mathbf{n}_1, \mathbf{n}_2, \ldots, \mathbf{n}_J\}$, such that an adding-up constraint $\sum_{i=1}^m n_{ij} = n$, for $j = 1, 2, \ldots, J$, is satisfied. $\mathbf{n}$ will be referred to as *occurrence vector*, $\mathcal{H}_n$ as *occurrence-vector working set*.

Let $\mathcal{P}$ be a set of all probability vectors, such that $\sum_{i=1}^m p_i = 1$.

Then, a simple question can be asked.

**Question 1.** *What is the most probable occurrence vector $\hat{\mathbf{n}}$, among occurrence vectors $\mathbf{n}$ from the working set $\mathcal{H}_n$, to be generated by the prior generator $\mathbf{q}$?*

The answer to the question is

$$(2.1) \qquad \hat{\mathbf{n}} = \arg\max_{\mathbf{n} \in \mathcal{H}_n} \pi(\mathbf{n}|\mathbf{q})$$

where

$$(2.2) \qquad \pi(\mathbf{n}|\mathbf{q}) = \frac{n!}{n_1! n_2! \ldots n_m!} \prod_{i=1}^m q_i^{n_i}$$

is the probability of generating the occurrence vector $\mathbf{n}$ by a prior generator $\mathbf{q}$.

The following theorem states the main result on asymptotic equivalence of Q1 and REM/MaxEnt. The proof is developed in the Appendix.

**Theorem 1.** *Let* $\mathbf{q}$ *be the prior generator and* $\mathcal{H}_n$ *be the working set. Let* $\hat{\mathbf{n}}$ *be the most probable occurrence vector from the working set* $\mathcal{H}_n$, *to be generated by the prior generator* $\mathbf{q}$. *And let* $n \to \infty$. *Then*

$$\frac{\hat{\mathbf{n}}}{n} = \hat{\mathbf{p}}$$

*where*

$$\hat{\mathbf{p}} = \arg\max_{\mathbf{p} \in \mathcal{P}} H(\mathbf{p}, \mathbf{q})$$

*and*

$$H(\mathbf{p}, \mathbf{q}) = -\sum_{i=1}^{m} p_i \ln\left(\frac{p_i}{q_i}\right)$$

*is the relative entropy of probability vector* $\mathbf{p}$ *on generator* $\mathbf{q}$.

**Corollary.** *If also some other* differentiable *constraint* $F(\frac{\mathbf{n}}{n}) = 0$ *is employed to form the working set* $\mathcal{H}_n$, *and a corresponding constraint* $F(\mathbf{p}) = 0$ *is added to the relative entropy maximization, the claim of Theorem 1 remains valid.*

*Note.* If the prior generator is uniform, $H(\mathbf{p}, \mathbf{q})$ reduces to Shannon's entropy $H(\mathbf{p}) = -\sum_{i=1}^{m} p_i \ln p_i$.

**Example.** Let $\mathbf{q} = [0.13\ 0.09\ 0.42\ 0.36]'$. Let $\mathcal{H}_n$ consists of all occurrence vectors $\{\mathbf{n}_1, \mathbf{n}_2, \ldots, \mathbf{n}_J\}$ such that

(2.3)
$$\sum_{i=1}^{m} n_{ij} = n \quad \text{for} \quad j = 1, 2, \ldots, J$$
$$\sum_{i=1}^{m} n_{ij} x_i = 3.2n \quad \text{for} \quad j = 1, 2, \ldots, J$$

where $\mathbf{x} = [1\ 2\ 3\ 4]'$.

Table 1 shows $\frac{\hat{\mathbf{n}}}{n}$ and $J$, for $n = 10, 50, 100, 500, 1000$, together with probability vector $\hat{\mathbf{p}}$ maximizing relative entropy $H(\mathbf{p}, \mathbf{q})$ under constraints

(2.4)
$$\sum_{i=1}^{m} p_i = 1$$
$$\sum_{i=1}^{m} p_i x_i = 3.2$$

Results for uniform prior generator are in the fourth column of the table.

Also, a different question can be asked.

**Question 2.** *What is the expected occurrence vector* $\bar{\mathbf{n}}$, *of occurrence vectors* $\mathbf{n}$ *from the working set* $\mathcal{H}_n$, *which can be generated by the prior generator* $\mathbf{q}$?

Specifying the notion of expectation,

(2.5)
$$\bar{\mathbf{n}} = \frac{\sum_{j=1}^{J} \pi(\mathbf{n}_j|\mathbf{q}) \mathbf{n}_j}{\sum_{j=1}^{J} \pi(\mathbf{n}_j|\mathbf{q})}$$

TABLE 1

| $n$ | $J$ | $\frac{\mathring{\mathbf{n}}}{n}|\mathbf{q}$ | $\frac{\mathring{\mathbf{n}}}{n}|$ uniform prior |
|---|---|---|---|
| 10 | 10 | 0.1000 0.0000 0.5000 0.4000 | 0.1000 0.1000 0.3000 0.5000 |
| 50 | 154 | 0.0800 0.0600 0.4400 0.4200 | 0.0800 0.1400 0.2800 0.5000 |
| 100 | 574 | 0.0800 0.0700 0.4200 0.4300 | 0.0800 0.1400 0.2800 0.5000 |
| 500 | 13534 | 0.0820 0.0700 0.4140 0.4340 | 0.0780 0.1460 0.2740 0.5020 |
| 1000 | 53734 | 0.0830 0.0700 0.4110 0.4360 | 0.0790 0.1460 0.2710 0.5040 |
| $\hat{\mathbf{p}}$ | | 0.0826 0.0709 0.4103 0.4361 | 0.0788 0.1462 0.2714 0.5037 |

**Theorem 2.** *Let* $\mathbf{q}$ *be the prior generator and* $\mathcal{H}_n$ *be the working set. Let* $\bar{\mathbf{n}}$ *be the expected occurrence vector of the working set* $\mathcal{H}_n$, *to be generated by the prior generator* $\mathbf{q}$. *And let* $n \to \infty$. *Then, using notation of Theorem 1,*

$$\frac{\bar{\mathbf{n}}}{n} = \hat{\mathbf{p}}$$

**Corollary.** *If also some other moment consistency constraint* $F(\frac{\mathbf{n}}{n}) = 0$ *is employed to form the working set* $\mathcal{H}_n$, *and a corresponding constraint* $F(\mathbf{p}) = 0$ *is added to the relative entropy maximization, the claim of Theorem 1 seems to remain valid – as following example indicates.*

*Note.* Scope of generality of the Theorem 2 is under both theoretical and numerical investigations. Numerical calculations show that the above Corollary is very likely to be valid.

**Example.** Assuming the same setup as in the previous Example, expected occurrence vectors for $n = [10, 50, 100, 1000]$ and generators $\mathbf{q} = [0.13\ 0.09\ 0.42\ 0.36]'$, $\mathbf{q} = [0.25\ 0.25\ 0.25\ 0.25]'$ are in Table 2.

TABLE 2

| $n$ | $\frac{\bar{\mathbf{n}}}{n}|\mathbf{q}$ | $\frac{\bar{\mathbf{n}}}{n}|$ uniform prior |
|---|---|---|
| 10 | 0.0721 0.0736 0.4365 0.4178 | 0.0701 0.1510 0.2877 0.4912 |
| 50 | 0.0806 0.0714 0.4153 0.4327 | 0.0771 0.1471 0.2745 0.5013 |
| 100 | 0.0816 0.0712 0.4128 0.4344 | 0.0779 0.1466 0.2729 0.5025 |
| 500 | 0.0824 0.0710 0.4108 0.4358 | 0.0786 0.1463 0.2717 0.5035 |
| 1000 | 0.0825 0.0709 0.4106 0.4360 | 0.0787 0.1462 0.2715 0.5036 |
| $\hat{\mathbf{p}}$ | 0.0826 0.0709 0.4103 0.4361 | 0.0788 0.1462 0.2714 0.5037 |

Note that $\frac{\bar{\mathbf{n}}}{n}$ converges to $\hat{\mathbf{p}}$ faster than $\frac{\mathring{\mathbf{n}}}{n}$.

*Comments.* Moving towards a summary of the Q1/Q2 rationale for $I$-divergence minimization we make following points:

1) Theorem 1 as any asymptotic theorem, can be interpreted in two directions: either moving towards infinity or moving back to finiteness. The first direction provides the rationale for REM. The second direction shows, that proper place of REM is only in the asymptotic.

2) Proof of Theorem 1 implies, that constraints that bind REM should be differentiable.

3) Hidden behind (2.1) and (2.2) is an assumption about *iid* sampling. Thus, REM as an asymptotic form of Q1, seems to be limited to the iid case.

4) Theorem 2, although its scope is yet investigated only numerically, suggests that REM is also asymptotic form of another elementary probabilistic question (Q2): *what is the expected vector of absolute frequencies (occurrences), among vectors allowed by a constraints, which can be generated by a prior generator?*.

### 3. Question 3 and Non-Parametric Maximum Likelihood

The second setup to be considered is following.

**Setup 2.** Let $\mathbf{n} = [n_1, n_2, \ldots, n_m]'$ be an occurrence vector.

Let $\mathcal{Q}$ be a set of prior generators. $\mathcal{Q}$ will be referred to as *prior-generator working set*.

Then, a simple question can be asked.

**Question 3.** *What is the most probable prior generator $\hat{\mathbf{q}}$, among all generators from the working set $\mathcal{Q}$, to generate the occurrence vector $\mathbf{n}$?*

Recalling (2) it is easy to see, that $\hat{\mathbf{q}}$, answering the question, solves

$$(3.1) \qquad \hat{\mathbf{q}} = \arg\max_{\mathbf{q} \in \mathcal{Q}} \prod_{i=1}^{m} q_i^{n_i}$$

**Example.** Let $\mathbf{n} = [13 \ 9 \ 42 \ 36]'$, and let $\mathcal{Q}$ consists of all generators $\mathbf{q}$ such that

$$\sum_{i=1}^{m} q_i = 1$$

$$\sum_{i=1}^{m} q_i x_i = 3.2$$

where $\mathbf{x} = [1 \ 2 \ 3 \ 4]'$.

The problem solved, for instance by means of Lagrangean, leads to $\hat{\mathbf{q}} = [0.0860 \ 0.0704 \ 0.4013 \ 0.4423]'$. Note that occurrence vectors $k\mathbf{n}$, for $k = 1, 2, \ldots$, produce the same $\hat{\mathbf{q}}$.

If $\mathbf{n}$ is taken a uniform one, then $\hat{\mathbf{q}} = [0.1009 \ 0.1384 \ 0.2204 \ 0.5403]'$.

The results are worth comparing with those of the previous Example.

*Comments.* Objective function of (3.1) can be called *non-parametric likelihood function*.

Note, that (3.1) is *identical* with maximization of relative entropy

$$\max_{\mathbf{q} \in \mathcal{Q}} H(\mathbf{p}, \mathbf{q})$$

Thus, non-parametric likelihood maximization is embedded in relative entropy maximization with respect to the prior generator. In this sense, relative entropy maximization with respect to generated pmf and relative entropy maximization with respect to generating pmf (identical with the non-parametric likelihood maximization) are 'inverse' tasks.

As it is well-known, parametric ML is also related to REM. In the parametric case REM on moment consistency constraints and ML on respective exponential form (family) lead to the same estimators of the pmf/pdf parameters (see [3]).

Question 1 and Question 3 can be repeated, in different order. But they can be asked also both at once. It is done in the third general setup.

## 4. Question 4 and Maximum Jeffreys' Entropy

**Setup 3.** Let $\mathcal{Q}$ be a prior-generator working set, $\mathcal{H}_n$ be an occurrence-vector working set. Let $\mathbf{q}$ be a prior generator from the generator-working set.

Then, a question might be asked.

**Question 4.** *What is the most probable occurrence vector $\tilde{\mathbf{n}}$, among all $\mathbf{n}$ vectors from the occurrence-vector working set $\mathcal{H}_n$,*

*1)* (the first setup) *to be generated by the prior generator* $\mathbf{q}$,

*and at the same time, (after a change of roles)*

*2)* (the second setup) *when taken as a generator $\frac{\mathbf{n}}{n}$, to generate occurrence vector* $\mathbf{q}n$?

Mathematically speaking,
$$\tilde{\mathbf{n}} = \arg\max_{\mathbf{n}\in\mathcal{H}_n} \pi\left((\mathbf{n}|\mathbf{q}), \left(\mathbf{q}n|\frac{\mathbf{n}}{n}\right)\right)$$

where
$$\pi\left((\mathbf{n}|\mathbf{q}), \left(\mathbf{q}n|\frac{\mathbf{n}}{n}\right)\right) = \pi(\mathbf{n}|\mathbf{q})\,\pi\left(\mathbf{q}n|\frac{\mathbf{n}}{n}\right)$$

By the same argument as in the case of Theorem 1 (see the Appendix), following Theorem can be proved.

**Theorem 3.** *Let $\mathbf{q}$ be the prior generator and $\mathcal{H}_n$ be the working set. Let*
$$\tilde{\mathbf{n}} = \arg\max_{\mathbf{n}\in\mathcal{H}_n} \pi\left((\mathbf{n}|\mathbf{q}), \left(\mathbf{q}n|\frac{\mathbf{n}}{n}\right)\right)$$

*and let $n \to \infty$. Then*
$$\frac{\tilde{\mathbf{n}}}{n} = \tilde{\mathbf{p}}$$

*where*
$$\tilde{\mathbf{p}} = \arg\max_{\mathbf{p}\in\mathcal{P}} J(\mathbf{p}, \mathbf{q})$$

*and*
$$J(\mathbf{p}, \mathbf{q}) = -\sum_{i=1}^{m} \left\{ p_i \ln\left(\frac{p_i}{q_i}\right) + q_i \ln\left(\frac{q_i}{p_i}\right) \right\}$$

*is Jeffreys' relative entropy of probability vector $\mathbf{p}$ on generator $\mathbf{q}$.*

The Jeffreys' relative entropy maximization reduces into
$$\max_{\mathbf{p}} \; EnLi(\mathbf{p}, \mathbf{q}) = \sum_{i=1}^{m} \left\{ q_i \ln p_i - p_i \ln\left(\frac{p_i}{q_i}\right) \right\}$$

For the important uniform prior generator case, $EnLi$ maximization further reduces into
$$\max_{\mathbf{p}} \; EnLi(\mathbf{p}) = \sum_{i=1}^{m} (1 - p_i) \ln p_i$$

**Example.** Let $\mathbf{q} = [0.13\ 0.09\ 0.42\ 0.36]'$. Let $\mathcal{H}_n$ consists of all occurrence vectors $\{\mathbf{n}_1, \mathbf{n}_2, \ldots, \mathbf{n}_J\}$ such that (2.3) is satisfied.

For $n = 10, 50, 100, 500, 1000$, $\tilde{\mathbf{n}}/n$ are stated in Table 3, together with the probability vector $\tilde{\mathbf{p}}$ maximizing Jeffreys' entropy under the constraints (2.4). Results for uniform prior generator are in the third column of the table.

TABLE 3

| $n$ | $\frac{\tilde{\mathbf{n}}}{n}\|\mathbf{q}$ | $\frac{\tilde{\mathbf{n}}}{n}\|$ uniform prior |
|---|---|---|
| 10 | 0.1000 0.1000 0.3000 0.5000 | 0.1000 0.1000 0.3000 0.5000 |
| 50 | 0.0800 0.0600 0.4400 0.4200 | 0.1000 0.1400 0.2200 0.5400 |
| 100 | 0.0800 0.0700 0.4200 0.4300 | 0.0900 0.1400 0.2500 0.5200 |
| 500 | 0.0840 0.0700 0.4080 0.4380 | 0.0920 0.1400 0.2440 0.5240 |
| 1000 | 0.0840 0.0710 0.4060 0.4390 | 0.0930 0.1390 0.2430 0.5250 |
| $\tilde{\mathbf{p}}$ | 0.0844 0.0705 0.4056 0.4394 | 0.0926 0.1395 0.2433 0.5246 |

Comparing $\hat{\mathbf{p}}$, $\hat{\mathbf{q}}$ and $\tilde{\mathbf{p}}$ shows, that $\tilde{\mathbf{p}}$ lays between $\hat{\mathbf{p}}$, $\hat{\mathbf{q}}$, what is implied by *EnLi*'s comprising of relative entropy and likelihood functions.

*Comments.* $J$-divergence (as is minus Jeffreys' entropy more commonly known) was introduced by [9] (see also [10]), mentioned in [12] and investigated further by [11]. $J$-divergence minimization (Jeffreys' entropy maximization) criterion seems to appear here, for the first time.

A question can be formulated in analogy to Question 2, in terms of expected occurrence vector

$$\tilde{\tilde{\mathbf{n}}} = \frac{\sum_{j=1}^{J} \pi\left((\mathbf{n}|\mathbf{q}), (\mathbf{q}n|\frac{\mathbf{n}}{n})\right) \mathbf{n}_j}{\sum_{j=1}^{J} \pi\left((\mathbf{n}|\mathbf{q}), (\mathbf{q}n|\frac{\mathbf{n}}{n})\right)}$$

and a theorem similar to the Theorem 2 can be stated.

## 5. SUMMARY

It was claimed, demonstrated and proved, that $I$-divergence minimization criterion is just an asymptotic form of simple probabilistic question, Question 1 of Section 2 (Theorem 1, Corollary, Example 1, Proof). Also, another probabilistic rationale, based on Question 2 was proposed and investigated, numerically (Theorem 2, Corollary, Example 2). Non-parametric maximum likelihood (Q3 method) was mentioned, and a discussion of relationship of both non-parametric and parametric likelihood methods to $I$-divergence minimization was presented. Finally, the fourth question was formulated, making the probabilistic ground for $J$-divergence minimization.

## Appendix A. Proof of Theorem 1

*Proof.* 1)

$$\max_{\mathbf{n}} \pi(\mathbf{n}|\mathbf{q})$$

subject to

$$\sum_{i=1}^{m} n_i = n$$

For the purpose of maximization $\pi(\mathbf{n}|\mathbf{q})$ can be *log*-transformed, into

$$v(\mathbf{n}) = \gamma(1 + \sum_{i=1}^{m} n_i) - \sum_{i=1}^{m} \gamma(n_i + 1) + l$$

where $\gamma(\cdot) = \ln \Gamma(\cdot)$, $\Gamma(\cdot)$ is gamma-function, and $l = \sum_{i=1}^{m} n_i \ln q_i$. Necessary condition for maximum of $\pi(\mathbf{n}|\mathbf{q})$ than is

(A.1) $$\mathrm{d}v(\mathbf{n}) = \sum_{i=1}^{m} [-\gamma'(n_i + 1) + \ln q_i] \, \mathrm{d}n_i = 0$$

since, according to the assumed adding-up constraint, $\sum_{i=1}^{m} \mathrm{d}n_i = 0$.

First, it will be proved that

(A.2) $$\lim_{n_i \to \infty} [\gamma'(n_i) - \ln(n_i + k)] = 0$$

for any $n_i$, and any $k$.

The first derivative of $\gamma$ can be written in a form of infinite series (see [2])

$$\gamma'(n_i) = g(n_i) - C$$

where

$$g(n_i) = \sum_{j=0}^{\infty} \left( \frac{1}{j+1} - \frac{1}{j+n_i} \right)$$

$$C = \text{Euler's constant}$$

For $n_i \in \mathbb{Z}$, denoted $\hat{n}_i$, the series $g(n_i)$ reduces into a harmonic series $H$

$$g(\hat{n}_i) = \sum_{j=1}^{\hat{n}_i - 1} \frac{1}{j} = H(\hat{n}_i - 1)$$

Let $n_i = \hat{n}_i + h$, where $0 \leq h < 1$. Then

$$\frac{1}{j+1} - \frac{1}{j+\hat{n}_i} \leq \frac{1}{j+1} - \frac{1}{j+\hat{n}_i+h} < \frac{1}{j+1} - \frac{1}{j+\hat{n}_i+1},$$

so

$$g(\hat{n}_i) = \sum_{j=1}^{\hat{n}_i - 1} \frac{1}{j} \leq g(n_i) < \sum_{j=1}^{\hat{n}_i} \frac{1}{j} = g(\hat{n}_i + 1).$$

Since, difference of the major series converges to zero,

$$\lim_{n_i \to \infty} [g(\hat{n}_i + 1) - g(\hat{n}_i)] = \lim_{n_i \to \infty} \frac{1}{\hat{n}_i} = 0$$

also

(A.3) $$\lim_{n_i \to \infty} [g(n_i) - g(\hat{n}_i)] = 0.$$

Due to a known property of harmonic series

$$\lim_{n_i \to \infty} [H(n_i) - \ln(n_i + k) - C] = 0, \text{ for any } k$$

holds also

(A.4) $$\lim_{n_i \to \infty} [g(\hat{n}_i) - \ln(n_i + k) - C] = 0$$

thus, adding (A.3), (A.4) up gives

$$\lim_{n_i \to \infty} [g(n_i) - \ln(n_i + k) - C] = 0,$$

respectively, for the derivative (recall that $\gamma'(n_i) = g(n_i) - C$)

$$\lim_{n_i \to \infty} [\gamma'(n_i) - \ln(n_i + k)] = 0,$$

what is just (A.2).

Without loss of generality, for $n \to \infty$ we can restrict for sub-space of $n_i \to \infty$, for $i = 1, 2, \ldots, m$, if on the sub-space exists the maximum, so the conditions of maximum for $n_i \to \infty$ take, thanks to (A.2), form

(A.5a) $$\sum_{i=1}^{m} [-\ln n_i + \ln q_i] \, dn_i = 0$$

(A.5b) $$\sum_{i=1}^{m} n_i = n$$

Due to Lemma 1, with $k = \frac{1}{n}$, system (A.5) can be transformed into an equivalent one

(A.6a) $$\sum_{i=1}^{m} \left[-\ln \frac{n_i}{n} + \ln q_i\right] d\frac{n_i}{n} = 0$$

(A.6b) $$\sum_{i=1}^{m} \frac{n_i}{n} = 1$$

2)

$$\max_{\mathbf{p}} H(\mathbf{p}, \mathbf{q})$$

subject to

$$\sum_{i=1}^{m} p_i = 1$$

Thus necessary conditions for maximum of relative entropy $H(\mathbf{p}, \mathbf{q})$, constrained by the respective adding-up constraint are

$$\text{(A.7a)} \qquad dH(\mathbf{p}, \mathbf{q}) = \sum_{i=1}^{m}[-\ln p_i + \ln q_i]dp_i = 0$$

$$\text{(A.7b)} \qquad \sum_{i=1}^{m} p_i = 1$$

Comparing (A.7) and (A.6) completes the proof. □

*Note.* Claim of the Corollary of Theorem 1 is immediately implied by the proof.

**Lemma 1.**
$$dv(k\mathbf{n}) = k\, dv(\mathbf{n})$$
*if* $\sum_{i=1}^{m} n_i = n$.

*Proof.* $dv(k\mathbf{n}) = \sum_{i=1}^{m}[-\ln kn_i + \ln q_i]dkn_i = k\sum_{i=1}^{m}[-\ln k - \ln n_i + \ln q_i]dn_i = k\sum_{i=1}^{m}[-\ln n_i + \ln q_i]dn_i = k\, dv(\mathbf{n})$. □

Institute of Measurement Science, Slovak Academy of Sciences, Dúbravská cesta 9, 842 19 Bratislava, Slovakia
*E-mail address*: umergren@savba.sk
*URL*: http://www.um.savba.sk/lab_15/mg.html

Railways of Slovak Republic, DDC, Klemensova 8, 813 61 Bratislava, Slovakia
*E-mail address*: grendar.marian@zsr.sk